\documentclass{jT}

\usepackage{amssymb,amsmath,latexsym,amsthm}
\usepackage[english]{babel}

\newcounter{theoremcounter}
\newcounter{lemmacounter}

\newcounter{dummycounter}

\newcounter{quescounter}
\newcounter{emptycounter}

\newtheorem{theorem}[theoremcounter]{Theorem}

\newtheorem{question}[quescounter]{Question}
\newtheorem{lemma}[lemmacounter]{Lemma}

\numberwithin{equation}{section}
\numberwithin{lemmacounter}{section}
\numberwithin{propcounter}{section}
\numberwithin{corcounter}{section}
\numberwithin{conjcounter}{section}
\numberwithin{theoremcounter}{section}
\numberwithin{defcounter}{section}
\numberwithin{probcounter}{section}

\newcounter{eqncounter}

\numberwithin{equation}{eqncounter}

\def\IF{\mathbb F}

\def\IZ{\mathbb Z}

\def\IP{\mathbb P}
\def\IQ{\mathbb Q}

\newenvironment{rproof}{\addvspace{\medskipamount}\par\noindent{\it Proof.\/}}
{\unskip\nobreak\hfill$\Box$\par\addvspace{\medskipamount}}
\renewcommand{\vec}[1]{\mbox{\boldmath$#1$}}
\def\ord{\mathop{\rm ord}\nolimits}
\def\Oseen{{\mathcal{O}}}
\def\A{{\mathfrak{A}}}

\def\c{{\mathfrak{c}}}

\def\B{{\mathfrak{B}}}

\def\L{L}
\def\k{k}
\def\K{K}
\def\M{F}
\def\vx{\vec {x}}
\def\MK{M(\K)}
\def\MM{M(\M)}

\def\v{\B}

\def\w{\wp}
\def\q{p}
\def\l{l}
\def\degK{\deg_{\K}}
\def\degM{\deg_{\M}}
\def\degk{\deg_{\k}}
\def\degR{\deg_{R}}
\def\d{d}
\def\kbar{\overline{\k}}

\def\comM{\hat{\M_{\w}}}

\begin{document}

\title[Small generators of function fields]{Small generators of function fields}


\author[Martin {\sc Widmer}]{{\sc Martin} Widmer}
\address{Martin {\sc Widmer}\\
Institut f\"ur Mathematik A\\
Technische Universit\"at Graz\\
Steyrergasse 30/II\\
8010 Graz\\ 
Austria}
\email{widmer@tugraz.at}

\maketitle

\begin{resume}
Soit $\K/\k$ une extension finie d'un corps global, donc $\K$ contient un \'el\'ement primitif $\alpha$, c'est \`a dire $\K=\k(\alpha)$.
Dans cet article, nous d\'emontrons l'existence d'un \'el\'ement
primitif de petite hauteur en cas d'un corps de fonctions. Notre r\'esultat r\'epond \`a une qu\'estion de Ruppert en cas d'un corps de fonctions.
\end{resume}

\begin{abstr}
Let $\K/\k$ be a finite extension of a global field. Such an extension can be generated over $\k$ by a 
single element. The aim of this article is to prove the existence of a ''small`` generator in the function field case.
This answers the function field version of a question of Ruppert on small generators of number fields.
\end{abstr}

\bigskip
\section{Introduction}

Let $\K$ be a finite extension of a global field $\k$ 
where global field means finite extension of either $\IQ$ or of a rational function field of transcendence degree one
over a finite field. 
Such an extension is generated by a single element and there exists a
natural concept of size on $\K$ given by the height. The well-known Theorem of Northcott
(originally proved for algebraic numbers but easily seen to hold also in positive characteristic)
implies that for each real $T$ there are only finitely many $\alpha \in \K$ whose height does not exceed $T$.
In particular there exists a smallest generator.
It is therefore natural to ask for lower and upper bounds for the height of a smallest generator.
We emphasize the situation where $d$ is fixed and $\K$ runs over all extensions of $\k$ of degree $d$.\\

Several people proved lower bounds for generators; first Mahler \cite{Mah} for the ground field $\k=\IQ$ and then
Silverman \cite{9} for arbitrary ground fields (and also higher dimensions), but see also
\cite{Rupp}, \cite{8} and \cite{EllVentorclass} for simpler results. For an extension $\K/\k$ of number fields Silverman's inequality 
implies 
\begin{alignat}1\label{Silverman}
h(1,\alpha)\geq \log|\Delta_{\K}|/(2d(d-1))-\log|\Delta_{\k}|/(2(d-1))-[\k:\IQ]\log d/(2(d-1))
\end{alignat}
for any generator $\alpha$ of $\K/\k$.
As shown by examples of Masser (Proposition 1 \cite{8}) and Ruppert \cite{Rupp}, 
this bound is sharp, at least up to an additive constant depending only on $\k$ and $d$.
A version of Silverman's bound in the function field case follows quickly from Castelnuovo's inequality.
For simplicity let us temporarily assume $\K$ and $\k$ are finite separable extensions of the rational function field $\IF_q(t)$ both with field of constants $\IF_q$.
We apply Castelnuovo's inequality as in \cite{Stichtenoth} III.10.3.Theorem with $F=\K=\k(\alpha)$, $F_1=\k$ and $F_2=\IF_q(\alpha)$.
Writing $g_{\k}$ and $g_{\K}$ for the genus of $\k$ and $\K$
we conclude $[\K:\IF_q(\alpha)]\geq g_{\K}/(d-1)-dg_{\k}/(d-1)+1$. From (\ref{degtransform}), (\ref{efd}) and the definition of the height in (\ref{height})
we easily deduce $h(1,\alpha)\geq [\K:\IF_q(\alpha)]/d$ and thus 
\begin{alignat}1\label{Castelnuovo}
h(1,\alpha)\geq g_{\K}/(d(d-1))-g_{\k}/(d-1)+1/d.
\end{alignat}
The discriminant $\Delta_{\K}=q^{\deg Diff(\K/\IF_q(t))}$ of $\K/\IF_q(t)$ is related with the genus by the Riemann-Hurwitz formula, more precisely  $\Delta_{\K}=q^{2g_{\K}+2([\K:\IF_q(t)]-1)}$. Thus (\ref{Castelnuovo}) matches with (\ref{Silverman}), at least up to an additive constant depending only on the degrees of $\k$ and $\K$.
A similar inequality as in (\ref{Castelnuovo}) was given by Thunder (\cite{Thslff} Lemma 6). Opposed to Silverman Thunder does not 
assume separability for the extension $\K/\k$.\\

What about upper bounds for the smallest generator? It seems that this problem has not been studied much yet.
However, at least for number fields the problem has been proposed explicitly by Ruppert. More precisely Ruppert 
(\cite{Rupp} Question 2) addressed the following question.
\begin{question}[Ruppert, 1998]
Does there exist a constant $C=C(d)$ such that for each number field $\K$ of degree $d$ there exists a generator $\alpha$
of the extension $\K/\IQ$ with $h(1,\alpha)\leq \log|\Delta_{\K}|/(2d)+C$?
\end{question}
Ruppert used the non-logarithmic naive height whereas we use the logarithmic projective absolute Weil height as defined in \cite{BG}.
However, it is easily seen that the question formulated here is equivalent to Ruppert's Question 2 in \cite{Rupp}.
One can show that there exists always an integral generator $\alpha$ of $\K/\IQ$ with $h(1,\alpha)\leq \log|\Delta_{\K}|/d$,
for a proof of this simple fact see \cite{VaalerWidmer}.
On the other hand, if $\alpha$ is an integral generator of an imaginary quadratic field $\K$ with minimal polynomial
$x^2+bx+c=(x-\alpha)(x-\overline{\alpha})$ then 
$h(1,\alpha)=\log\sqrt{\alpha\overline{\alpha}}=\log \sqrt{c}\geq \log\sqrt{4c-b^2}-\log 2\geq \log|\Delta_{\K}|/d-\log 2$.
Nevertheless, Ruppert showed that Question 1 has an affirmative answer for $\k=\IQ$ and $\K$ either quadratic or a totally real field of 
prime degree. In fact, using Minkowski's convex body Theorem to construct a Pisot-number generator, it suffices to assume $\K$ has a real embedding and one can drop the prime degree condition. For more details we refer
to \cite{VaalerWidmer}. Ruppert's result for $d=2$ relies heavily on a distribution result of Duke \cite{Duke} that does not appear to have an analogue 
for higher degrees and is ineffective. As a consequence Ruppert's constant $C$ for $d=2$ is ineffective.\\

In this note we introduce a completely different strategy which applies in the function field case and the number field case.
However, it relies on the existence of a certain divisor which is guaranteed under GRH but might be rather troublesome
to establish unconditionally.
The aim of this short note is to answer positively Ruppert's question in the function field case.
So let $\k$ be an algebraic function field with finite constant field $\k_0$ and transcendence degree one over $\k_0$. We have the following result.
\begin{theorem}\label{thffc}
Let $\K$ be a finite field extension of $\k$.
There exists an element $\alpha$ in $\K$ with $\K=\k(\alpha)$ and 
a constant $C=C(\k,[\K:\k])$
depending solely on $\k,[\K:\k]$ such that
\begin{alignat*}1
h(1,\alpha)\leq \frac{g_K}{\d(\K/\k)}+C
\end{alignat*}
where $g_{\K}$ denotes the genus of the function field $\K$ with field of constants $\K_0$ and $\d(\K/\k)=[\K:\k]/[\K_0:\k_0]$.
\end{theorem}

\section{Notation and definitions}
Throughout this note we fix an algebraic closure $\kbar$ of $\k$. All fields are considered to be subfields of $\kbar$.
For any finite extension $\M$ of $\k$ we write $\M_0$ for the field of constants in $\M$;
in other words $\M_0$ is the algebraic closure of $\k_0$ in $\M$.
When we talk of the field $\M$ we implicitly mean the field $\M$ with field of constants $\M_0$.
We define the geometric degree $\d(\M/\k)$ of the extension
$\M$ over $\k$ as
\begin{alignat*}1
\d(\M/\k)=\frac{[\M:\k]}{[\M_0:\k_0]}.
\end{alignat*}
Let $\MM$ be the set of all places in $\M$. For a place $\w$ in $\MM$ let
$\Oseen_{\w}$ be the valuation ring of $\M$ at $\w$; we can identify $\w$ with the 
unique maximal ideal in $\Oseen_{\w}$. 
We write $\M_{\w}=\Oseen_{\w}/\w$ for the residue class field and 
$\comM$ for the topological completion of $\M$ at the place $\w$.
Write $\ord_{\w}$ for the order function on $\comM$ normalized to have image in $\IZ\cup \infty$.
We extend $\ord_{\w}$ to $\comM^n$ by defining
\begin{alignat*}1
\ord_{\w}(x_1,...,x_n)=\min_{1\leq i\leq n}\ord_{\w}x_i
\end{alignat*}
with the usual convention $\ord_{\w}0=\infty>0$.
Each non-zero element $\vx$ of $\M^n$ gives rise to a divisor $(\vx)$ over $\M$
\begin{alignat*}1
(\vx)=\sum_{\w\in \MM}\ord_{\w}(\vx)\cdot \w.
\end{alignat*}
For a divisor $A$ over $\M$ we define the Riemann-Roch space in $\M^n$
\begin{alignat*}1
\L_n(A)=\{\vx\in \M^n\backslash 0;(\vx)+A\geq 0\}\cup \{0\}.
\end{alignat*}
This is a $\M_0$
vector space of finite dimension. Denote its dimension over $\M_0$ by $\l_n(A)$.\\

The degree of a place $\w$ in $\MM$ is defined by $\degM \w=[\M_{\w}:\M_0]$. Let $\K$ be a finite
extension of $\M$ and let $\v$ be a place in $\MK$ above $\w$. We write 
$f(\v/\w)=[\K_{\v}:\M_{\w}]$ for the residue degree of $\v$ over $\w$. Then we have  
\begin{alignat}1\label{degtransform}
\nonumber \degK \v=[\K_{\v}:\K_0]=\frac{[\K_{\v}:\M_0]}{[\K_0:\M_0]}&=\frac{[\K_{\v}:\M_{\w}]}{[\K_0:\M_0]}[\M_{\w}:\M_0]\\
&=\frac{f(\v/\w)}{[\K_0:\M_0]}\degM \w.
\end{alignat}
Writing $e(\v/\w)$ for the ramification index we also have 
\begin{alignat}1\label{efd}
\sum_{\v\mid \w}e(\v/\w)f(\v/\w)=[\K:\M],
\end{alignat}
see for example III.1.11.Theorem in \cite{Stichtenoth}.\\

Each divisor $A=\sum_{\w}a_{\w}\w$ over the smaller field $\M$ naturally defines a divisor
\begin{alignat*}1
A^{(\K)}=\sum_{\w}\sum_{\v\mid \w}a_{\w}e(\v/\w)\v
\end{alignat*}
over the large field $\K$.\\

As in \cite{Thslff} we define the height $h$ on non-zero $\vx$ in $\K^n$ by 
\begin{alignat}1
\label{height}
h(\vx)=-\frac{\degK (\vx)}{\d(\K/\k)}.
\end{alignat}
Note that the degree of a principal divisor is zero so that the height defines a function on
projective space $\IP^{n-1}(\K)$ over $\K$ of dimension $n-1$. This shows also that the height is nonnegative
since to evaluate the height of $\vx$ we can assume that one coordinate is $1$.
Moreover it is absolute in the following sense. Suppose $\vx \in \K^n$ and let $D$ be the divisor given by $D=(\vx)$.
Let $R$ be a finite extension of $\K$ and view $\vx \in R^n$. Let $D^{(R)}$ be the divisor over $R$ given by $D^{(R)}=(\vx)$.
By \cite{Artin} Chap.15, Thm.9 we have $\degR(D^{(R)})=\d(R/\K)\degK(D)$ and
by \cite{Artin} Chap.15, Thm.2 we have $\d(R/\k)=\d(R/\K)\d(\K/\k)$. Thus $h(\vx)$ remains unchanged if one views
$\vx$ in $R^n$. Therefore the height extends to a projective height on $\kbar^n$.
Suppose $\vx \in \L_n(\A)\subseteq \K^n$ then directly from the definition we see that 
\begin{alignat}1
\label{heightineq}
h(\vx)\leq \frac{\degK \A}{\d(\K/\k)}.
\end{alignat}

\section{The strategy}\label{strat}
Let $S$ be a finite set of places in $M(\K)$ such that the following two properties hold:

\noindent(i) for each place $\q$ in $M(\k)$ there is at most one place in $S$ that lies above $\q$,\\
(ii)$f(\v/\q)=1$ for all $\v\in S$ and $\q\in M(\k)$ with $\v\mid \q$.\\ 

A set $S$ with these two properties will be called admissible. 
Note that for each field $F$ with $\k\subseteq \M\subseteq \K$
and for all places $\v,\v'\in S$, $\w,\w'\in M(\M)$ with $\v\mid \w$ and $\v'\mid \w'$ we have
\begin{alignat}1
\label{uniqueext}
&\v\neq \v' \Rightarrow \w\neq \w',\\
\label{degreeone}
&f(\v'/\w')=1.
\end{alignat}
We say the divisor $\A$ is admissible if it can be written in the form
\begin{alignat}1\label{suitdiv}
\A=\sum_{\v\in S}1\cdot\v.
\end{alignat}
with an admissible set $S$.
\begin{lemma}\label{lemma1}
Suppose $\A$ is an admissible divisor and suppose $\vx=(1,x)$ with $x \notin \K_0$ and $\vx \in \L_2(\A)$. Then $\k(x)=\K$ and $h(\vx)\leq \degK \A/\d(\K/\k)$.
\end{lemma}
\begin{rproof}
Suppose $\k(x)=\M \subsetneq \K$ and write $(\vx)=\sum_{\w}a_{\w}\w$ for the divisor over $\M$. 
When we consider $(\vx)$ as a divisor over $K$ we have
$(\vx)=\sum_{\v}a_{\v}\v$ with $a_{\v}=a_{\w}e(\v/\w)$. 
Note that none of the coefficients $a_{\w}$ is positive and since $x\notin \K_0$ at least one is negative, say $a_{\w'}$. 
Since $(\vx)$ lies in $\L_2(\A)$ and $\A$ is an admissible divisor we conclude by (\ref{uniqueext}) that
there is exactly one place $\v'$ in $M(\K)$ with $\v'\mid\w'$ and by (\ref{degreeone}) that $f(\v'/\w')=1$. 
Together with (\ref{efd}) we deduce that
$1<[\K:\M]=\sum_{\v'\mid \w'}e(\v'/\w')f(\v'/\w')=e(\v'/\w')$. This means that $a_{\v'}=a_{\w'}e(\v'/\w')<a_{\w'}\leq-1$, contradicting the fact $\A+(\vx)\geq 0$.
Thus $\k(x)=\K$. The remaining statement comes from (\ref{heightineq}).
\end{rproof}
\begin{lemma}\label{lemma2}
Suppose $\A$ is an admissible divisor with $\degK \A>g_{\K}$.
Then there exists $\alpha$ in $\K$ with $\k(\alpha)=\K$ and 
\begin{alignat*}1
h(1,\alpha)\leq \frac{\degK \A}{\d(\K/\k)}.
\end{alignat*}
\end{lemma}
\begin{rproof}
We apply the Theorem of Riemann-Roch to the space $\L_1(\A)=\{x\in \K\backslash 0;(x)+\A\geq 0\}\cup \{0\}$ 
to conclude $\l_1(\A)>1$. 
Therefore we find a $\alpha$ in $\L_1(\A)\backslash \K_0$.
Now $(\vx)=(1,\alpha)$ and $(\alpha)$ have the same pole-divisors and since $\A\geq 0$
we see that $(\vx)$ lies in $\L_2(\A)$.
Applying Lemma \ref{lemma1} proves the lemma.
\end{rproof}

\section{Constructing a suitable divisor}
In this section we will prove that $\K$ admits an admissible divisor of degree $g_{\K}+1$
provided $g_{\K}$ is large enough.
\begin{lemma}\label{lemma3}
Let $l$ be a positive integer. The number of places $\v\in M(\K)$ with
\begin{alignat*}1
\degK \v&=l \text{ and }\\
f(\v/\q)&=1 \text{ for } \q\in M(\k) \text{ with } \v\mid \q
\end{alignat*}
is
\begin{alignat*}1
\geq \frac{|\K_0|^{l}}{l}-(2+7g_{\K})\frac{|\K_0|^{l/2}}{l}-l[\K_0:\k_0][\K:\k](|\K_0|^{l/2}+(2+7g_{\k})|\K_0|^{l/4}).
\end{alignat*}
\end{lemma}
\begin{rproof}
Using Riemann's hypothesis one can obtain a good lower bound for the 
number of places $\v$ of fixed degree. For instance V.2.10 Corollary (a) in \cite{Stichtenoth}
tells us that the total number of places $\v\in M(\K)$ with $\degK \v=l$ is 
\begin{alignat*}1
\geq \frac{|\K_0|^{l}}{l}-(2+7g_{\K})\frac{|\K_0|^{l/2}}{l}.
\end{alignat*}
From (\ref{degtransform}) we get
$\degk \q=l[\K_0:\k_0]/f(\v/\q)$ for $\q\in M(\k)$, $\v\mid \q$.
Suppose $f=f(\v/\q)>1$. 
Applying V.2.10 Corollary (a) again we get the following upper bound for the number 
of places $\q\in M(\k)$ with $\degk \q=l[\K_0:\k_0]/f$
\begin{alignat*}1
&\frac{f|\k_0|^{l[\K_0:\k_0]/f}}{l[\K_0:\k_0]}+(2+7g_{\k})\frac{f|\k_0|^{l[\K_0:\k_0]/(2f)}}{l[\K_0:\k_0]}\\
=&\frac{f|\K_0|^{l/f}}{l[\K_0:\k_0]}+(2+7g_{\k})\frac{f|\K_0|^{l/(2f)}}{l[\K_0:\k_0]}\\
\leq& |\K_0|^{l/2}+(2+7g_{\k})|\K_0|^{l/4}.
\end{alignat*}
Above each place $\q$ in $M(\k)$ there lie at most $[\K:\k]$ places $\v\in M(\K)$.
Therefore the number of places $\v\in M(\K)$ satisfying $\degK \v=l$ and $f(\v/\q)=f>1$ is
\begin{alignat*}1
\leq [\K:\k](|\K_0|^{l/2}+(2+7g_{\k})|\K_0|^{l/4}).
\end{alignat*}
Summing over all divisors $f$ of $l[\K_0:\k_0]$ we find that the
number of places $\v\in M(\K)$ satisfying $\degK \v=l$ and $f(\v/\q)>1$ is
\begin{alignat*}1
\leq l[\K_0:\k_0][\K:\k](|\K_0|^{l/2}+(2+7g_{\k})|\K_0|^{l/4}).
\end{alignat*}
This in turn implies that the number of places $\v\in M(\K)$ satisfying $\degK \v=l$ and $f(\v/\q)=1$ is 
\begin{alignat*}1
\geq \frac{|\K_0|^{l}}{l}-(2+7g_{\K})\frac{|\K_0|^{l/2}}{l}-l[\K_0:\k_0][\K:\k](|\K_0|^{l/2}+(2+7g_{\k})|\K_0|^{l/4})
\end{alignat*}
and this proves the lemma.
\end{rproof}

\section{Proof of Theorem \ref{thffc}}
Let $d$ be a positive integer. 
We can find a constant $C_1=C_1(\k,d)$ such that with $l=g_{\K}+1$ 
\begin{alignat*}1
\frac{|\K_0|^{l}}{l}-(2+7g_{\K})\frac{|\K_0|^{l/2}}{l}-l[\K_0:\k_0][\K:\k](|\K_0|^{l/2}+(2+7g_{\k})|\K_0|^{l/4})>0
\end{alignat*}
for all extensions $\K$ of $\k$ satisfying $g_{\K}>C_1$ and $[\K:\k]=d$.
By virtue of Lemma \ref{lemma3} we conclude that for each such $\K$
there exists a place
$\v\in M(\K)$ with $\degK\v=g_{\K}+1$ and $f(\v/\q)=1 \text{ for } \q\in M(\k) \text{ with } \v\mid \q$.
In particular there exists an admissible divisor, namely $\v$, with $\degK\v=g_{\K}+1$.
From Lemma \ref{lemma2} we conclude that if $g_{\K}>C_1$ and $[\K:\k]=d$ then there exists $\alpha \in \K$ with $\K=\k(\alpha)$
and $h(1,\alpha)\leq (g_{\K}+1)/\d(\K/\k)$. There are only finitely many field extensions $\K$ of $\k$  of degree $d$ with $g_{\K}\leq C_1$. 
Hence there exists a constant $C$ as in Theorem \ref{thffc}
depending solely on $C_1,\k,d$ and thus depending solely on $\k,d$ 
such that the statement of Theorem \ref{thffc} holds for all extensions $\K$ of $\k$ of degree $d$.

\section*{Acknowledgements}
I would like to thank Felipe Voloch for giving me the idea of the proof of Lemma \ref{lemma3} and pointing out the 
connection with Castelnuovo's bound mentioned in the Introduction.
Moreover I am grateful to David Masser, Jeff Thunder and Jeff Vaaler for many delightful discussions.

\bibliographystyle{amsplain}
\bibliography{literature}

\end{document}